\documentclass[11pt]{article}
\usepackage{latexsym}
\usepackage{amsmath}
\usepackage{amssymb}
\newtheorem{theorem}{Theorem}[section]

\newtheorem{lemma}[theorem]{Lemma}

\def\p{{\bf Proof.} \quad}
\def\q{\hfill\rule{1ex}{1ex}}

\begin{document}
\title{\bf On two $q$-ary $n$-cube coloring problems}
\author{{{\small\bf Zhe
Han}\thanks{email: hanzhe101@gmail.com.}\quad{\small\bf Mei Lu}\thanks{email: mlu@math.tsinghua.edu.cn.}}\\
{\small Department of Mathematical Sciences, Tsinghua
University, Beijing 100084, China}}

\date{}

\maketitle\baselineskip 16.3pt

\begin{abstract}
Let $\chi'_d(n,q)$ (resp. $\chi_d(n,q)$) denote the minimum number of colors necessary to color a $q$-ary $n$-cube
so that no two vertices that are at a distance at most $d$ (resp. exactly $d$) get the same color. These two problems were proposed in the study of scalability of optical networks. In this paper, we provide upper and lower bounds on $\chi'_d(n,q)$ and $\chi_d(n,q)$ when $q$ is a prime power.

\end{abstract}

\vskip.3cm

\section{Introduction}

Let $q=p^m$, where $p$ is a prime.
Let $V_n$ be the $n$-dimensional vector space over a finite field $\mathbb{F}_q$, i.e.,
$$V_n=\{(x_1,\ldots,x_n)| x_i\in \mathbb{F}_q\}.$$
\noindent For $\mathbf{x}=(x_1,\ldots,x_n)\in V_n$ and
$\mathbf{y}=(y_1,\ldots,y_n)\in V_n$, the Hamming distance $d_H(\mathbf{x},\mathbf{y})$ between
$\mathbf{x}$ and $\mathbf{y}$ is the number of the coordinates in which they differ, i.e.,
$$d_H(\mathbf{x},\mathbf{y})=|\{i|x_i\neq y_i,~1\le i\le n\}|.$$
\noindent The Hamming weight $w_H(\mathbf{x})$ of a vector $\mathbf{x}\in V_n$ is the number
of nonzero coordinates in $\mathbf{x}$. Obviously,
$$d_H(\mathbf{x},\mathbf{y})=w_H(\mathbf{x}-\mathbf{y}).$$

Let $G=(V ,E)$ be a simple undirected graph with the vertex set $V$ and the edge set $E$. For two distinct vertices $\mathbf{u}$ and $\mathbf{v}$ in $V$, the distance between $\mathbf{u}$ and $\mathbf{v}$, denoted by $d(\mathbf{u},\mathbf{v})$, is the number of edges in a shortest path joining $\mathbf{u}$ and $\mathbf{v}$ and the diameter of $G$ is
the maximum distance between any two vertices of $G$.

 A {\em $q$-ary $n$-cube} \cite{Bhuyan}, denoted by $Q_n(q)$, is an undirected graph with the vertex set $V_n$ and the edge set
 $$E_n=\{(\mathbf{x},\mathbf{y})|\mathbf{x},\mathbf{y}\in V_n,d_{H}(\mathbf{x},\mathbf{y})=1\}.$$
 The 2-ary $n$-cube $Q_n(2)$ is a hypercube suggested by Sullivan and Bashkow \cite{Sullivan} and is one of the most popular, versatile and efficient topological structures of interconnection networks. By the definition, $Q_n(q)$ is $n(q-1)$-regular of order $q^n$, $d(\mathbf{x},\mathbf{y})=d_H(\mathbf{x},\mathbf{y})$ for any distinct vertices $\mathbf{x},\mathbf{y}\in V_n$ and the diameter of $Q_n(q)$ is $n$. Note that $Q_n(2)$ is a bipartite, but $Q_n(q)$ is not bipartite when $q\ge 3$.

A coloring of $V_n$ with $L$ colors is a map $\Gamma$ from the vertex set $V_n$ to ${\cal L}=\{1,2,\ldots,L\}$.
A $d$-distance (resp. exactly $d$-distance) coloring of $V_n$ is to color the vertices of $V_n$ such that
any two vertices with Hamming distance at most $d$ (resp. exactly $d$) have different colors. Note that
for a coloring of $V_n$ with $L$ colors
\[\Gamma: V_n \longrightarrow \mathcal{L},\]
 it is a $d$-distance coloring of $V_n$ if and only if for any two distinct vertices $\mathbf{x},\mathbf{y}\in V_n$,
\[\Gamma(\mathbf{x})\neq \Gamma(\mathbf{y}) \mbox{ if } d_H(\mathbf{x},\mathbf{y})\leq d\]
\noindent and it is an exactly $d$-distance coloring of $V_n$ if and only if for any two distinct vertices $\mathbf{x},\mathbf{y}\in V_n$,
\[\Gamma(\mathbf{x})\neq \Gamma(\mathbf{y}) \mbox{ if } d_H(\mathbf{x},\mathbf{y})= d.\]
Denote $\chi'_d(n,q)$ (resp. $\chi_d(n,q)$) as the minimum number of colors needed for a $d$-distance (resp. an exactly
$d$-distance) coloring of $V_n$. Clearly, $\chi_{d}(n,q)\leq \chi'_{d}(n,q)$ for each $n$.

These two coloring problems originally arose in the study of the scalability of optical networks \cite{Pavan}. Some bounds on $\chi'_d(n,2)$ and $\chi_d(n,2)$ were given and some exact values were determined (see for example \cite{Enomoto}-\cite{Ostergard}, \cite{Skupie}-\cite{Wan}). In the paper, we will give some lower and upper bounds on $\chi_{d}(n,q)$ and $\chi'_{d}(n,q)$. In particular, some bounds are given by using methods in coding theory.




The $d$-distance coloring and exactly $d$-distance coloring of $V_n$ are equivalent to certain partitions of $V_n$,
which are related
to codes in coding theory. So we first introduce some definitions in coding theory. A nonempty subset $C$ of $V_n$ is called a $q$-ary code of length $n$. Any element in $C$ is called a codeword of $C$. The minimum distance $d(C)$ is defined as the minimum Hamming distance between two distinct codewords of $C$. A $q$-ary code $C$ of length $n$ and minimum distance $d$ is called an $(n, d)_q$ code. A $q$-ary code $C$ of length $n$ and minimum distance at least $d$ is called an $(n,\ge d)_q$ code. Let $A_q(n,d)$ (resp. $A_q(n,\ge d)$) denote the maximum size of an $(n, d)_q$ (resp. $(n,\ge d)_q$) code. A $q$-ary code $C$ is called a $q$-ary $[n, k]_q$ linear code if $C$ is a $k$-dimensional subspace of $V_n$. A $q$-ary $[n,k]_q$ linear code with minimum distance $d$ is called a $q$-ary $[n,k,d]_q$. For the fixed $n$ and minimum distance $d$, let $k(n,d)_q$ denote the maximum dimension of a $q$-ary $[n,k,d]_q$ code. A $q$-ary code $C$ of length $n$ is called an $(n, \overline{\{d\}})_q$ forbidden distance code if $d_H(\mathbf{u},\mathbf{v})\not= d$ for any two distinct codewords $\mathbf{u},\mathbf{v}\in C$. Given $n$ and $d$, let $Q(n, d)_q$ denote the maximum size of a $q$-ary $(n, \overline{\{d\}})_q$ forbidden distance code.

An $(n,L,d)_q$-partition (resp. $(n,L,\overline{\{d\}})_q$-partition) of $V_n$ is a set of subsets $\{B_i\}_{i=1}^{L}$ of $V_n$ satisfying (i) $B_i\bigcap B_j=\emptyset$ for $i\neq j$ and $\bigcup_{i=1}^{L}B_i=V_n$; (ii) each $B_i$ is a $q$-ary $(n,\geq d)_q$ code (resp. $(n,\overline{\{d\}})_q$ forbidden distance code). It is well known that a $d$-distance coloring (resp. an exactly $d$-distance) of $V_n$ with $L$ colors is equivalent to an $(n,L,d+1)_q$-partition (resp. $(n,L,\overline{\{d\}})_q$-partition) of $V_n$. Hence $\chi'_{d}(n,q)$ (resp. $\chi_{d}(n,q)$) is the minimum number $L$ of subsets in any $(n,L, d+1)_q$ (resp. $(n,L,\overline{\{d\}})_q$-partition) of $V_n$.

Since $q^n=|V_n|=\sum_{i=1}^{L}|B_i|\le L A_{q}(n,\geq d+1)\leq L A_{q}(n,d+1) $ and $A_q(n,d)$ is decreasing in $d$, we have
\begin{equation}\label{eq:basic}
 \chi'_{d}(n,q)\geq \frac{q^n}{A_{q}(n,d+1)}.
 \end{equation}

Note that for a $q$-ary $[n,k,d+1]_q$ linear code $C$, the cosets of $C$ form an $(n,q^{n-k},d+1)_q$-partition of $V_n$. Hence,
if there exists a $q$-ary $[n,k,d+1]_q$ linear code, then
\begin{equation} \label{eq:code}
\chi'_d(n,q)\leq q^{n-k}.
 \end{equation}
In particular,
\begin{equation} \label{eq:max}
\chi'_d(n)\leq q^{n-k(n,d+1)}.
 \end{equation}
Furthermore, if $A_q(n,d+1)=q^{k(n,d+1)}$, i.e., $A_q(n,d+1)$ is attained by a $q$-ary $[n,k,d+1]_q$ linear code, then
\begin{equation}\label{eq:final}
\chi'_{d}(n,q)=q^{n-k(n,d+1)}.
\end{equation} Since $q$-ary Hamming code $[\frac{q^r-1}{q-1},\frac{q^r-1}{q-1}-r,3]_q$ and $q$-ary simplex code $[\frac{q^r-1}{q-1},r,q^{r-1}]_q$ are optimal code, by (4) we have that
\begin{equation} \chi'_{2}\left(\frac{q^r-1}{q-1},q\right)=q^r,\end{equation}and
\begin{equation} \chi'_{q^{r-1}-1}\left(\frac{q^r-1}{q-1},q\right)= q^{\frac{q^r-1}{q-1}-r}. \end{equation}

In the next two Sections, we give some lower and upper bounds on $\chi_{d}(n,q)$ and $\chi'_{d}(n,q)$ for general $d$.


\section{Bounds of $\chi'_{d}(n,q)$ }

First, we present the lower bound of $\chi'_{d}(n,q)$.

\begin{theorem} $\chi'_{d}(n,q)\ge \sum\limits_ {i=0}^{d/2}(q-1)^i{n\choose i}$ if $d$ is even and $\chi'_{d}(n,q)\ge q\sum\limits_{i=0}^{(d-1)/2}(q-1)^i{n-1\choose i}$ if $d$ is odd.
\end{theorem}

\p  Suppose $d$ is even. We consider all vertices within distance $d/2$ from the vertex $(0,0,\cdots,0)$. The total number of these vertices is $\sum\limits_ {i=0}^{d/2}(q-1)^i{n\choose i}$. The distance
between any two of them is at most $d$. Thus, these vertices should have distinct colors.

Now we suppose $d$ is odd. we consider all vertices within distance $(d-1)/2$ of the vertex  $(a,0,\cdots,0)$ for any $a\in \mathbb{F}_q$. The total number of these vertices
is $q\sum\limits_{i=0}^{(d-1)/2}(q-1)^i{n-1\choose i}$. The distance between any two of them is at most $d$. Thus, these vertices
should have distinct colors.
\q

\vskip.2cm

{\bf Note.} By Theorem 2.1, $\chi'_{d}(n,2)\ge \sum\limits_ {i=0}^{d/2}{n\choose i}$ for even $d$ and $\chi'_{d}(n,2)\ge 2\sum\limits_{i=0}^{(d-1)/2}{n-1\choose i}$ for odd $k$, which had been shown by Kim et al in \cite{Kim}. In addition,
Theorem 2.1 can be proved by Singleton bound in coding theory easily.

\vskip.2cm

Next, we give some upper bounds of $\chi'_{d}(n,q)$.

\vskip.2cm

\begin{theorem}For integers $n$ and $d$, we have

\begin{equation}
 \chi'_{d}(n,q) \leq     q^{\left\lfloor\log_q\left(\sum_{j=0}^{d-1}(q-1)^{j}{n-1\choose j}\right)\right\rfloor+1}.
\end{equation}
\end{theorem}

\p First, we show that there exists a $t\times n$ matrix $H$ with the property that any $d$ columns of $H$ are linearly independent over $\mathbb{F}_q$ and there exist $d+1$ dependent columns. This $H$ is actually the parity check matrix of an $[n,n-t,d+1]_q$-code.

Now, we describe a procedure for constructing a $t\times n$ parity check matrix $H$ by choosing its column vectors sequentially.
The first column vector can be any nonzero vector. Suppose that we already have a set $V$ of $i$ vectors so that
any $d$ of them are linearly independent. The $(i+1)$th vector can be chosen as long as it is not the linear combination of any
$d-1$ vectors in $V$. In other words, since we are working over the field $\mathbb{F}_q$, the total number of unallowable vectors is at most $(q-1){i\choose 1}+(q-1)^2{i\choose 2}+\ldots+(q-1)^{d-1}{i\choose d-1}$. Consequently, as long as $(q-1){i\choose 1}+(q-1)^2{i\choose 2}+\ldots+(q-1)^{d-1}{i\choose d-1}<q^{t}-1$, we can still add a new column to $H$.
Therefore a $t\times n$ parity check matrix $H$ can be constructed.  It follows that
$t=\left\lfloor\log_q\left(\sum_{j=0}^{d-1}(q-1)^{j}{n-1\choose j}\right)\right\rfloor+1$ and so $\chi'_{d}(n,q)\leq q^{t}$ by (\ref{eq:code}). \q
\vskip.2cm

\vskip.2cm

Now we consider the bounds of $\chi'_{d}(n,q)$ for $d=2$.

For a positive integer $k$ and a prime $l$, a $k$-dimensional $l$-adic representation of a nonnegative integer $i$ is a  $k$-dimensional vector
$(x_1, x_2,\cdots, x_k)$, where $x_i\in \{0,1,\ldots,l-1\}$,  such that $$x_1\cdot q^{k-1} + x_2\cdot q^{k-2}+\cdots + x_{k-1}\cdot q+x_k = i.$$

\vskip.2cm

\begin{theorem} $(q-1)n+1\le \chi'_{2}(n,q)$ and $\chi'_{2}(n,q)\le q^{\lceil\log_qn\rceil+1}$ if $q=p$.
\end{theorem}

\p Consider vertices $(0,0,\cdots,0)$, $(a,0,\cdots,0)$, $(0,a,\cdots,0)$, $\cdots$, and $(0,\cdots,0,a)$, where $a\in \mathbb{F}_q\setminus \{0\}$. Then the distance between any two vertices is at most $2$. Therefore, they should have distinct colors.
That is, $\chi'_{2}(n,q)\ge (q-1)n+1$.

Let $m=\lceil\log_qn\rceil+1$ and $M_q(n)$ be an $n\times m$ matrix whose $(i + 1)$th row $r_i$ ($0\le i\le n-1$) has the form $(b_i, n_i)$, where $b_i$ is the $(m-1)$-dimensional $q$-adic representation $(x_{i1}, x_{i2},\cdots, x_{i,m-1})$ of integer $i$ and $n_i=\sum_{j=1}^{m-1}x_{ij}+1(\mod q)$. For example,
$$M_3(5)=\left(
  \begin{array}{ccc}
    0 & 0 & 1 \\
    0 & 1 & 2 \\
    0 & 2 & 0 \\
    1 & 0 & 2 \\
    1 & 1 & 0 \\
  \end{array}
\right).$$
Denote colors by $m$-dimensional vectors over $\mathbb{F}_q$. For any vertex $\mathbf{v}=(v_1,v_2,\cdots,v_n)$, color $\mathbf{v}$ with the color $\mathbf{v}M_q(n)$, where all arithmetic operations are performed in $\mathbb{F}_q$. Then there are  $q^m$ colors in total. Let $\mathbf{u}$ and $\mathbf{v}$ be any two vertices with distance $d\le 2$, and use $r_i=(x_{i1}, x_{i2},\cdots, x_{i,m-1},n_i)$  to denote the $(i+1)$th row of $M_q(n)$, where $0\le i\le n-1$. We consider two cases.

{\bf Case 1} $d=1$.

Suppose $u_i\not=v_i$ for some $1\le i\le n$. Then $\mathbf{u}M_q(n)-\mathbf{v}M_q(n)=(u_i-v_i)r_{i-1}$. By the definition of $M_q(n)$, $r_{i-1}\not=\mathbf{0}$. Thus $\mathbf{u}M_q(n)\not=\mathbf{v}M_q(n)$.

{\bf Case 2} $d=2$.

Suppose $u_i\not=v_i$ and $u_j\not=v_j$ for some $1\le i<j\le n$. Then $$\mathbf{u}M_q(n)-\mathbf{v}M_q(n)=(u_i-v_i)r_{i-1}+(u_j-v_j)r_{j-1}.$$ Denote $u_i-v_i=a_i$ and $u_j-v_j=a_j$ for short. Note that $(x_{i1}, x_{i2},\cdots, x_{i,m-1})\not=(x_{j1}, x_{j2},\cdots, x_{j,m-1})$ for $i\not=j$. If $a_i(x_{i1}, x_{i2},\cdots, x_{i,m-1})=a_j(x_{j1}, x_{j2},\cdots, x_{j,m-1})$, then $a_i\not=a_j$ but $a_i\sum_{t=1}^{m-1}x_{it}=a_j\sum_{t=1}^{m-1}x_{jt}$. Since $n_i=\sum_{j=1}^{m-1}x_{ij}+1(\mod q)$, we have $a_in_i\not=a_jn_j$, which implies $\mathbf{u}M_q(n)\not=\mathbf{v}M_q(n)$.

Therefore, if the distance of $\mathbf{u}$ and $\mathbf{v}$ is at most 2, then $\mathbf{u}M_q(n)\not=\mathbf{v}M_q(n)$, and so vertices $\mathbf{u}$ and $\mathbf{v}$ have different
colors. Hence $\chi'_{2}(n,q)\le q^{\lceil\log_qn\rceil+1}$.\q

\vskip.2cm

{\bf Note.} In \cite{Wan}, Wan showed that $\chi'_{2}(n,2)\le 2^{\lceil\log_2(n+1)\rceil}$. But the inequality  $\chi'_{2}(n,q)\le q^{\lceil\log_q(n+1)\rceil}$ does not always hold for $q$-ary $n$-cube. For example,  consider the case $q=3$ and $n=2$. Then $3^{\lceil\log_3(2+1)\rceil}=3$. But $\chi'_{2}(2,3)= 3^2>3$ since the distance of any two vertices in $Q_2(3)$ is at most 2 and so they should have distinct colors.

\section{Bounds of $\chi_{d}(n,q)$ }

In this section, we consider the bounds of $\chi_{d}(n,q)$. Since the diameter of $Q_n(q)$ is $n$, $\chi_{d}(n,q)= 1$ for any $d\ge n+1$. So we can assume $d\le n$.

\vskip.2cm

{

\begin{theorem} $\chi_{d}(n,q)\ge q$. Particularly, $\chi_{n}(n,q)= q$ and $\chi_{1}(n,q)= q$.
\end{theorem}

\p Consider vertices $(x_1,x_2,\cdots,x_n)$ with $x_1=x_2=\cdots=x_d=a$ for any $a\in \mathbb{F}_q$ and $x_{d+1}=\cdots=x_n=0$. Then the distance between any two vertices is $d$. Therefore, they should have distinct colors.
That is, $\chi_{d}(n,q)\ge q$.

For any $a\in\mathbb{F}_q$, denote $S_{a}=\{(x_1,x_2,\cdots,x_n)~:~x_1=a,x_j\in \mathbb{F}_q\mbox{~for~}2\le j\le n\}$. Then $V_n=\bigcup_{a\in \mathbb{F}_q}S_a$.  For vertices $\mathbf{u},\mathbf{v}\in S_a$ and some $a\in \mathbb{F}_q$, the distance between them is at most $n-1$ and they could have the same color. Thus $\chi_{n}(n,q)\le q$, which implies $\chi_{n}(n,q)= q$.

We consider the case when $d=1$. Let $\mathbb{F}_q=\{0, \alpha^0, \ldots,\alpha^{q-2}\}$ where $\alpha$ is a primitive element of $\mathbb{F}_q$ and ${\cal L}=\{0,1,\ldots,q-1\}$. Put $\Gamma(0,0,\ldots,0)=0$,

$$\Gamma(x_1,\ldots,x_n)=\sum_{x_k\neq 0, x_k=\alpha^{i_k},1\leq k\leq n}i_k+1\mod q,$$
where $x_1,\ldots,x_n\in \mathbb{F}_q$.
\noindent For any vertices $\mathbf{x}$ and $\mathbf{y}$ with distance $1$,  they have the different colors by the above definition. Thus $\chi_{1}(n,q)\le q$, which implies $\chi_{1}(n,q)= q$.  \q

{\bf Note.} When $d<n$ and $d>1$, $\chi_{d}(n,q)=q$ does not always hold since we can give two counterexamples. When $d=3, n=5, q=3$, we can find four vertices $(0,0,0,0,0)$,$(1,1,1,0,0)$, $(2,2,2,0,0)$, and $(2,0,1,2,0)$ whose distance between any two vertices is $3$, which means that these four vertices would have $4$ different colors. Likewise, when $d=5, n=7, q=3$, there are four vertices $(0,0,0,0,0,0,0)$, $(1,0,1,0,1,1,1)$, $(0,0,2,2,2,2,1)$, and $(0,1,1,1,2,1,0)$ whose distance between any two vertices is $3$, which means that these four vertices would have $4$ different colors.
But the conclusion holds for $q=2$ since $Q_n(2)$ is bipartite.

\vskip.2cm

\begin{theorem} For even $d$, $\chi_{d}(n,q)\ge \max\{q,\lfloor 2n/d\rfloor\}$.
\end{theorem}
\p Consider vertices $(\underbrace{a,\cdots,a}_{d/2},0,\cdots,0)$, $(\underbrace{0,\cdots,0}_{d/2},\underbrace{a,\cdots,a}_{d/2},0,\cdots,0)$, $(\underbrace{0,\cdots,0}_{d},$ $\underbrace{a,\cdots,a}_{d/2},0,\cdots,0)$$\cdots$, where $a\in \mathbb{F}_q\setminus\{0\}$. Then the distance between any two vertices is $d$. Therefore, they should have distinct colors.\q

\vskip.2cm

 In the following we give some results about the exactly $d$-distance coloring of $V_n$ by using the
 argument of linear forbidden distance codes.  A $q$-ary $[n,k]_q$ linear code $C$ is called a $q$-ary
$[n,k,\overline{d}]_q$ linear forbidden distance code if $C$ is also an $(n, \overline{d})$ forbidden
distance code, i.e.,  $w_{H}(\mathbf{c})\neq d$ for any nonzero codeword $\mathbf{c}$ of $C$.
Given $n$ and $d$, the maximum dimension of a $q$-ary
$[n,k,\overline{d}]_q$ linear forbidden distance code is denoted by $k(n,\overline{d})$. Similar to $d$-distance
coloring, we have

\begin{equation}
\chi_{d}(n,q)\leq q^{n-k}.
\end{equation}
\noindent
\begin{equation}\label{eq:exactlc}
\chi_{d}(n,q)\leq q^{n-k(n,\overline{d})}.
\end{equation}

In the following lemma the lower bound on $k(n,\overline{d})$ is given.

\vskip.2cm
\begin{lemma} We have
\begin{equation}
k(n,\overline{d})\geq n-\left\lceil \log_{q}[2+{n-1\choose d-1}(q-1)^{d-1}]\right\rceil.
\end{equation}
\end{lemma}
\p
First if
\begin{equation}
2+{n-1 \choose d-1}(q-1)^{d-1} \leq q^{m},
\end{equation}
then there exists a $q$-ary matrix $H$ of size $m\times n$ such that any column of $H$ is not equal to
the $\mathbb{F}_q$-linear combination of any other $d-1$ columns of $H$.
Hence there is a $q$-ary $[n,k,\overline{d}]_q$ linear forbidden distance code, where the dimension $k\geq n-m$. Take
\[m=\left\lceil \log_{q}[2+{n-1\choose d-1}(q-1)^{d-1}]\right\rceil. \]
Thus,
\begin{eqnarray*}
k(n,\overline{d}) &\geq & dim(C)\\
         & \geq & n-m \\
         & \geq & n-\left\lceil \log_{q}[2+{n-1\choose d-1}(q-1)^{d-1}]\right\rceil.
\end{eqnarray*}

This completes the proof.
\q

\vskip.2cm

\begin{theorem} We have
\begin{equation}
\chi_{d }(n,q)\leq q^{\lceil \log_{q}[2+{n-1\choose d-1}(q-1)^{d-1}]\rceil}.
\end{equation}
\end{theorem}

\p
The theorem follows from Lemma 3.3 and (\ref{eq:exactlc}).
\q



\vskip.2cm

\section*{Acknowledgements} This work is partially supported by National Natural Science Foundation of China
(Nos. 61373019, 11171097 and 61170289).

\end{document}